\documentclass[12pt,a4paper,oneside]{amsart}
\usepackage{amsfonts, amsmath, amssymb, amsthm, amscd}

\usepackage{anysize}

\usepackage[english]{babel}
\usepackage{graphicx}
\usepackage{url}
\usepackage{caption}
\usepackage{float}
\usepackage{tikz-cd}
\usepackage{wrapfig}

\usepackage{paralist} 
\usepackage{enumitem}

\usepackage{hyperref}
\hypersetup{
   colorlinks   = true, 
   urlcolor     = blue, 
   linkcolor    = blue, 
   citecolor   = red 
}

\def\Z{{\mathbb Z}}
\def\R{{\mathbb R}}

\newtheorem{thm}{Theorem}[section]

\newtheorem{clm}[thm]{Claim}
\newtheorem{cor}[thm]{Corollary}
\newtheorem{que}[thm]{Question}

\newtheorem*{thm*}{Theorem}

\theoremstyle{remark}

\newtheorem{rem}[thm]{Remark}

\theoremstyle{definition}
\newtheorem{define}[thm]{Definition}

\numberwithin{equation}{section}

\title{A subexponential size triangulation of $\R P^n$}

\author{Karim~Adiprasito{$^\spadesuit$}}

\author{Sergey~Avvakumov{$^\clubsuit$}}

\author{Roman~Karasev{$^\diamondsuit$}}

\thanks{{$^\spadesuit$} Has received funding from the European Research Council under the European Union’s Seventh Framework Programme ERC Grant agreement ERC StG 716424 - CASe and the Israel Science Founda-tion under ISF Grant 1050/16.}

\thanks{{$^\clubsuit$} Has received funding from the European Research Council under the European Union’s Seventh Framework Programme ERC Grant agreement ERC StG 716424 - CASe.}

\thanks{{$^\diamondsuit$} Supported by the Federal professorship program grant 1.456.2016/1.4 and the Russian Foundation for Basic Research grants 18-01-00036 and 19-01-00169}

\address{{Karim Adiprasito}, Einstein Institute of Mathematics, Hebrew University of Jerusalem, Jerusalem, Israel \emph{and} Department of Mathematical Sciences, University of Copenhagen, Copenhagen, Denmark}
\email{ka@math.ku.dk, adiprasito@math.huji.ac.il}

\address{Sergey~Avvakumov, Department of Mathematical Sciences, University of Copenhagen, Copenhagen, Denmark}
\email{savvakumov@gmail.com}

\address{Roman~Karasev, Moscow Institute of Physics and Technology, Institutskiy per. 9, Dolgoprudny, Russia 141700\newline \indent
Institute for Information Transmission Problems RAS, Bolshoy Karetny per. 19, Moscow, Russia 127994}
\email{r\_n\_karasev@mail.ru}
\urladdr{http://www.rkarasev.ru/en/}

\begin{document}

\maketitle

\begin{center}
\emph{Dedicated to Wolfgang K\"uhnel on the occasion of his 70th birthday.}
\end{center}
\begin{abstract}
We break the exponential barrier for triangulations of real projective space, constructing a trianglation of $\R P^n$ with $e^{(\frac{1}{2}+o(1))\sqrt{n}{\log n}}$ vertices.
\end{abstract}

\section{Introduction}

While in general, every smooth manifold allows for a triangulation, that is, a simplicial complex homeomorphic to the manifold in question, it is a natural and often hard problem to construct small triangulations of manifolds, and usually poses a difficult challenge. And so, outside of special cases, there are few cases known where the upper bounds and lower bounds come even close to each other.  

Let us focus on minimality in terms of the number of vertices that a triangulation of a given manifold would have, a study Banchoff and K\"uhnel initiated \cite{Kuhnel}. The best lower bounds in this area are usually either homological, or homotopic in nature. Indeed, 
\begin{compactitem}[$\circ$]
\item it is clear that the number of vertices cannot be lower than the ball-category, or the more studied Lusternik-Schnirelmann category \cite{CLOT}. In particular
\item it is bounded from below in terms of the cup length of the space in question. In fact, it is easy to show and observed by Arnoux and Marin that for a space of cup length $n$, one needs $\binom{n+2}{2}$ vertices \cite{AM}. Finally
\item Murai gave a lower bound in terms of the Betti numbers of (closed and orientable) manifolds \cite{MT}, which was simplified and generalized to general manifolds by Adiprasito and Yashfe \cite{AHL, AY}. This bound in general is not so good for interesting manifolds, as it seems insensitive to any interesting multiplicative structure in the cohomology ring, let alone homotopy.
\end{compactitem}

Returning to the example of $\R P^n$: Out of these available bounds, the one by Arnoux and Marin is best, yielding a lower bound of $\binom{n+2}{2}$ vertices \cite{AM}. Still, the best construction so far is essentially K\"uhnel's observation that the barycentric subdivision of the $n$-simplex yields a triangulation of the $(n-1)$-sphere on $2^{n+1}-2$ vertices, with a $\Z_2$-action, such that antipodal vertices are at distance at least $3$ from each other. This yields a triangulation of $\R P^{n-1}$ of size $2^{n}-1$.

Since then, no substantial improvement has been made for the general problem, and focus has shifted to experimental study of low-dimensional cases (see \cite{LutzFH} for an excellent survey) and improvements of the base of the exponential \cite{BasakSarkar, VZ}.

Surprisingly perhaps, and at least counter to the prevailing expectation of experts, we are not only constructing a triangulated sphere, but a polytope, counter to intuition coming from concentration inequalities on the sphere, that often imply exponential lower bounds \cite{Barvinok}.

\subsection*{Acknowledgments} We wish to thank Wolfgang K\"uhnel and Alexis Marin for very extensive and helpful remarks on the first version of the paper and anonymous referees for their comments and their suggestion of how to simplify the last part of the proof of Claim \ref{clm:sat}. We also wish to thank Arseniy Akopyan, Xavier Goaoc,  Andrey Kupavskii, Gaku Liu, J{\'a}nos Pach and D{\"o}m{\"o}t{\"o}r P{\'a}lv{\"o}lgyi for helpful discussions.

\section{Main result}

\begin{thm}
\label{thm:main}
For all positive integers $n$, there exists a convex centrally symmetric $n$-dimensional polytope $P$ such that:
\begin{itemize}
\item All the vertices of $P$ lie on the unit sphere and belong to either the positive coordinate orthant (with all coordinates non-negative) or the negative coordinate orthant (with all coordinates non-positive).
\item For any vertex $A\in P$, if $F, F'\subset P$ are facets with $A\in F$ and $-A\in F'$, then $F\cap F'=\emptyset$.
\item The number of vertices of $P$ is less than $e^{(\frac{1}{2}+o(1))\sqrt{n}{\log n}}$.
\end{itemize}
\end{thm}

Only the second property is required for the subsequent quotient of $\partial P$ being a \emph{triangulation} of $\R P^{n-1}$, as opposed to being just homeomorphic to $\R P^{n-1}$; otherwise a centrally symmetric $n$-polytope with only $2n$ vertices exists -- the crosspolytope.

Our proof indeed relies on starting with the crosspolytope, but then finding a smaller triangulation than the barycentric subdivision to achieve the desired property, exploiting that the barycentric subdivision is quite wasteful in adding vertices. Instead, we note that several vertices of the barycentric subdivision can be removed, without affecting the crucial second property. The boundary $\partial P$ can be viewed as a Delaunay triangulation \cite{Delone} of the sphere on the set of vertices of $P$ with the distance measured by the inner product of the corresponding unit vectors. This observation inspired parts of our proof of Theorem \ref{thm:main}, for example see property (3) in Claim \ref{clm:req}.

\begin{cor}
\label{thm:cor}
For all positive integers $n$, there exists a triangulation of $\R P^{n-1}$ with at most $e^{(\frac{1}{2}+o(1))\sqrt{n}{\log n}}$ vertices.
\end{cor}
\begin{proof}
Let $P$ be the polytope whose existence is guaranteed by Theorem \ref{thm:main}. The group $\Z_2$ acts on $\partial P$ by central symmetries. The polytope might not be simplicial, so before passing to the quotient $\partial P/\Z_2$ we need to triangulate its boundary first.

Let us construct a $\Z_2$-equivariant, i.e., symmetric, triangulation $S$ of $\partial P$. One way to do so is to take a so-called \emph{pulling} triangulation obtained by consequently pulling pairs of opposite vertices. Equivalently, shift the vertices of $P$ slightly so that any subset of them not containing opposite vertices is in general position. The movement can be chosen small enough so that the convex hull of the shifted vertices is still a triangulation of $\partial P$. The triangulated sphere $S$ has the same number of vertices as $P$.

From the second property of $P$ we get that the closed stars of any two opposite vertices of $S$ are disjoint. Hence, the quotient $S/\Z_2$ is a simplicial complex.
This complex is homeomorphic to $\R P^{n-1}$.
\end{proof}

\section{Proof of Theorem \ref{thm:main}}

Let $V$ be a subset of the set of non-empty subsets of $\{1,\dots,n\}$.
We identify any $A\in V$ with a unit vector in $\R^n$ whose endpoint has its $i$th coordinate equal to $1/{\sqrt{|A|}}$ if $i\in A$ and $0$ otherwise.

A centrally symmetric polytope $P(V)$ is the convex hull of the endpoints of the vectors $V\sqcup -V$.

In the following three claims we present an explicit construction of such $V$ that $P(V)$ satisfies the conditions of Theorem \ref{thm:main} thus proving the theorem. By $\langle \cdot,\cdot\rangle$ we denote the inner product.

\begin{clm}
\label{clm:req}
Suppose that $V$ satisfies the following properties:
\begin{itemize}
\item[(1)] $\{i\}\in V$ for all $i\in\{1,\dots,n\}$.
\item[(2)] If $A\in V$ and $|A|>1$, then $A\setminus i\in V$ for any $i\in A$.
\item[(3)] Let $A, B\in V$ and $X\in S^{n-1}$ be unit vectors such that $\langle A,B\rangle=0$ and $\langle A,X\rangle=\langle B,X\rangle>0$. Then there is $C\in V$ such that $\langle C,X\rangle>\langle A,X\rangle=\langle B,X\rangle$.
\end{itemize}

Then $P(V)$ satisfies the conditions of Theorem \ref{thm:main}, except maybe the condition on the number of vertices.
\end{clm}

\begin{proof}
By the property (1), $P(V)$ is an $n$-dimensional polytope.

Let $A\in V$ be a vertex and $F, F'\subset P(V)$ be facets with $A\in F$ and $-A\in F'$. We need to prove that $F\cap F'=\emptyset$.
Assume to the contrary, that there is a vertex $B\in F\cap F'$. Without loss of generality, we may assume that $B\in V$ (as opposed to $B\in-V$).

Let us check that by the property (2) for any vertex $D\in V\sqcup -V$ of $P(V)$ and any coordinate hyperplane $H \not\ni D$ the projection of $D$ to $H$ lies in the interior of $P(V)\cap H$. Because of the central symmetry of $P(V)$ we may assume that $D\in V$.
Hyperplane $H$ corresponds to the $i$th coordinate for some $i\in\{1,\dots,n\}$. Then $D \not\in H$ means that the $i$th coordinate of $D$ is non-zero (and even positive, because $D\in V$) when $D$ is considered as a unit vector; or equivalently that $i\in D$ when $D$ is considered as a subset of $\{1,\dots,n\}$. If $D=\{i\}$ then it projects to the origin. If the cardinality of $D$ is greater than $1$, then consider $D':=D\setminus i\in V$. The $i$th coordinate of the vector $D'$ is $0$ and the rest of its coordinates are equal to the corresponding coordinates of $D$ uniformly scaled by a factor $>1$ because $D'$ is unit. So, $D$ projects into the interior of the line segment connecting the origin with the endpoint of $D'$.

Suppose there is a coordinate hyperplane $H$ such that $B$ and $-A$ do not belong to $H$ and lie on the opposite sides of $H$.
Then both the projections of $B$ and $-A$ lie in the interior of $P(V)\cap H$.
So, the intersection point of $H$ with the line segment connecting $B$ with $-A$ also lies in the interior of $P(V)\cap H$.
On the other hand, the intersection point lies in the facet $F'\subset \partial P(V)$, contradiction.

So, there is no coordinate hyperplane strictly separating $B$ from $-A$.
On the other hand, all coordinates of $B$ are non-negative and all coordinates of $-A$ are non-positive.
Which means that $\langle B,-A\rangle=0$ and so $\langle A,B\rangle=0$.

Let $X$ be the outer normal of the facet $F\ni A,B$, clearly $\langle A,X\rangle=\langle B,X\rangle>0$. We already know that $\langle A,B\rangle=0$. So, we can use the property (3) to find $C\in V$ such that $\langle C,X\rangle>\langle A,X\rangle=\langle B,X\rangle$. The existence of such vertex $C$ contradicts to $F$ being a facet of $P(V)$.
\end{proof}

\begin{define}
\label{def:tight}
A subset $V$ of the set of non-empty subsets of $\{1,\dots,n\}$ is called \emph{tight} if

\begin{itemize}
\item[(1)] $\{i\}\in V$ for all $i\in\{1,\dots,n\}$.
\item[(2)] If $A\in V$ and $|A|>1$, then $A\setminus i\in V$ for any $i\in A$.
\item[(3)] For every $A,B\in V$ with $A\cap B=\emptyset$, there are $i\in A$ and $j\in B$ such that either 
\begin{itemize}
\item[(3a)] $B\sqcup i\in V$ and $A\sqcup j\setminus i\in V$,
\end{itemize}
or
\begin{itemize}
\item[(3b)] $A\sqcup j\in V$ and $B\sqcup i\setminus j\in V$.
\end{itemize}
\end{itemize}
 
\end{define}

\begin{clm}
\label{clm:sat}
Any tight subset $V$ satisfies the requirements of Claim \ref{clm:req}.
\end{clm}

\begin{proof}
Requirements (1) and (2) of Claim \ref{clm:req} are the same as the properties (1) and (2) in the definition of tightness.
So, we only need to check the requirement (3) of Claim \ref{clm:req}.

Let $A, B\in V$ and $X\in S^{n-1}$ be such that $\langle A,B\rangle=0$ and $\langle A,X\rangle=\langle B,X\rangle>0$.
We need to find $C\in V$ such that $\langle C,X\rangle>\langle A,X\rangle=\langle B,X\rangle$.

Denote the coordinates of $X$ by $x_1,\dots,x_{n}$.
Without loss of generality we may assume that $A=\{1,\dots,a\}$ and that $x_1\leq x_2\leq \dots \leq x_a$
(Indeed, we can find a permutation $\sigma$ such that $A=\{\sigma(1),\dots,\sigma(a)\}$ and then additionally ask from $\sigma$ that $x_{\sigma(1)}\leq x_{\sigma(2)}\leq \dots \leq x_{\sigma(a)}$. Then just omit $\sigma$ from the notation.). Note, that we cannot assume $x_1 \geq 0$.

From $\langle A,B\rangle=0$ we get that $A\cap B=\emptyset$.
So, without loss of generality we may assume that $A$ and $B$ satisfy the property (3b) in the definition of tightness, i.e.,  there are $i\in A$ and $j\in B$ such that $A\sqcup j\in V$ and $B\sqcup i\setminus j\in V$.

If $x_i > x_j$, then $C:=B\sqcup i\setminus j\in V$ is as required because $\langle C,X\rangle > \langle B,X\rangle$.
So, we can assume that $x_j \geq x_i \geq x_1$.

Denote $A_1:=A\setminus 1$ and $A_2:=A\sqcup j\in V$.

Assume first that $a=1$ and hence $A_1=\emptyset\notin V$. Then $\langle X, A \rangle=x_1$ meaning that $x_1 > 0$. So, 
$\langle X, A_2 \rangle=\frac{x_1+x_j}{\sqrt{2}}\geq\frac{2x_1}{\sqrt{2}}>x_1=\langle X, A \rangle$ and $C:=A_2$ is as required.

From now on we can assume that $a>1$ and $A_1\in V$.

We have $$\langle X, A \rangle = \frac{x_1+\dots+x_a}{\sqrt{a}},$$ $$\langle X, A_1 \rangle = \frac{x_2+\dots+x_a}{\sqrt{a-1}},$$ $$\langle X, A_2 \rangle =\frac{(x_1+x_j)+x_2+\dots+x_a}{\sqrt{a+1}} \geq \frac{2x_1+x_2+\dots+x_a}{\sqrt{a+1}}.$$

From the following inequalities it follows that either $C:=A_1$ or $C:=A_2$ is as required:

\[
\langle X, A \rangle =\frac{x_1+\dots+x_a}{\sqrt{a}} < \frac{x_1+\dots+x_a}{\frac{1}{2}(\sqrt{a-1}+\sqrt{a+1})} \leq {\rm max}(\langle X, A_1\rangle, \langle X, A_2\rangle),
\]
where
\begin{itemize}
\item the first inequality holds since $\frac{1}{2}(\sqrt{a-1}+\sqrt{a+1})<\sqrt{a}$ and $x_1+\dots+x_a> 0$ because $\langle X, A \rangle>0$,
\item and the second inequality holds because $\frac{\alpha+\gamma}{\beta+\delta}\leq {\rm max}(\frac{\alpha}{\beta}, \frac{\gamma}{\delta})$ for any reals $\alpha,\beta,\gamma,\delta$ with $\beta,\delta > 0$.
\end{itemize}

\end{proof}

\begin{clm}
\label{clm:siz}
There is a tight subset $V$ of size at most $e^{(\frac{1}{2}+o(1))\sqrt{n}{\log n}}$.
\end{clm}

\begin{rem}
After the first preprint of this paper was published, Peter Frankl, J{\'a}nos Pach, and D{\"o}m{\"o}t{\"o}r P{\'a}lv{\"o}lgyi proved that there are no tight subsets of size less than $2^{\sqrt{2n}}$ (in fact, they used some weaker assumptions than tightness), see \cite[Theorem 1]{frankl2021exchange} for details. In this sense the construction provided here is close to optimal.
\end{rem}

\begin{proof}
Partition the set $\{1,\dots,n\}$ into several disjoint subsets $S_1,\dots,S_k$.
Let $V$ be the set of non-empty subsets of $\{1,\dots,n\}$, whose intersection with each $S_i$, except maybe one of them, contains not more than one element.
I.e., $A\in V$ if and only if $A\neq\emptyset$ and there exists at most one $i\in \{1,\dots,k\}$ such that $|A\cap S_i|>1$.

Obviously, $V$ satisfies the properties (1) and (2) in the definition of tightness. Let us check that $V$ also satisfies (3a) or (3b), or sometimes both.

Let $A, B\in V$ be such that $A\cap B=\emptyset$. Let $S_a\in\{S_1,\dots,S_k\}$ be one of the subsets which maximizes $|A\cap S_a|$. Note, that the union of $A$ with any element of $S_a$ is still in $V$. Likewise, let $S_b$ be one of the subsets which maximizes $|B\cap S_b|$. Again, the union of $B$ with any element of $S_b$ is still in $V$.

Suppose that either $A\cap S_b$ or $B\cap S_a$ is non-empty. Without loss of generality, let us assume that $A\cap S_b$ is non-empty. Pick any elements $i\in A\cap S_b$ and $j\in B\cap S_b$. Then $i$ and $j$ satisfy (3a).

The remaining case is when both $A\cap S_b$ and $B\cap S_a$ are empty. Pick any elements $i\in A\cap S_a$ and $j\in B\cap S_b$. Then $i$ and $j$ satisfy both (3a) and (3b).

It remains to choose the subsets $S_1,\dots,S_k$ so that the size of $V$ is as required. We choose them to be roughly equal in size, each subset not larger than $s:=\lceil \frac{n}{k} \rceil$.

Let us analyze the size of $V$. For an element $A\in V$ we have
\begin{itemize}
\item $k$ choices of $S_a\in\{S_1,\dots,S_k\}$ which will maximize $|A\cap S_a|$,
\item at most $2^s$ choices of which elements of $S_a$ to add to $A$,
\item at most $s+1$ choices for each of $k-1$ subsets $S_j\neq S_a$ -- we can either add any one of at most $s$ elements of $S_j$ to $A$ or add none of them.
\end{itemize}

Note, that this way we count some elements of $V$ multiple times. In total we get
$$|V| < k\cdot 2^s\cdot(s+1)^{k-1}.$$

Choosing $k=\lceil \sqrt{n}\rceil$, we get the required asymptotic

$$|V| < (\sqrt{n}+1)\cdot 2^{\sqrt{n}}\cdot(\sqrt{n}+1)^{\sqrt{n}} = e^{(\frac{1}{2}+o(1))\sqrt{n}{\log n}}.$$

\end{proof}

\begin{rem}
Our choice of $k\approx s\approx \sqrt{n}$ in the proof of Claim \ref{clm:siz} is slightly suboptimal. Consequently, the bound on the number of vertices in Theorem \ref{thm:main} can be slightly improved, but only by a subpolynomial factor, which to state precisely we deemed not relevant.
\end{rem}

It should be noted that a subset does not have to be tight to satisfy the requirements of Claim \ref{clm:req}. 

In the end we would like to pose two questions which we do not know the answers to:

\begin{que}
Is there a non-tight subset which satisfies the requirements of Claim \ref{clm:req} and which is significantly smaller, i.e., of order at most $e^{n^{(\frac{1}{2}-O(1))}}$, than the tight subset we construct in Claim \ref{clm:siz}?
\end{que}

\begin{que}
Is there a better than quadratic in $n$ lower bound on the number of vertices of a convex centrally symmetric $n$-dimensional polytope $P$ satisfying the second property in the statement of Theorem \ref{thm:main}? What about both the first and the second property?
\end{que}

	{\small
		\bibliographystyle{myamsalpha}
		\bibliography{ref}}


\end{document}